\magnification 1050
 \def\oui{oui}
 \ifx\textures\oui
 \input CrayolaColors 
 \long\def\rge#1{\Red#1\Black}

 \else 
\def\arXiv{oui}
  
 \def\oui{oui} 
  
\ifx\arXiv\oui
\else
 \pdfpagewidth=210truemm
 \pdfpageheight=297truemm 
\fi
  
  %
  %

  %
  %

  \catcode`@=12 

 \def\defrefnote#1{\definexref{#1}{{\the\footnotenumber}}{refnotes}}

  %
  %


\ifx\couleurs\oui
\input graphicx
\fi

\input eplain.tex
 
\ifx\couleurs\oui
\beginpackages
\usepackage{color}
\endpackages 
\long\def\rge#1{{\color{red}#1}}

\definecolor{bleu-iecn}{cmyk}{.98,.13,.1,.55}

\fi

\makeatletter
\def\numberedfootnote{%
ÊÊ\global\advance\footnotenumber by 1
ÊÊ\@eplainfootnote{{\number\footnotenumber}}%
}%
\def\makecolumns#1/#2 {\par \begingroup
ÊÊ \@columndepth = #1
ÊÊ \advance\@columndepth by -1
ÊÊ \divide \@columndepth by #2
ÊÊ \advance\@columndepth by 1
ÊÊ \@linestogoincolumn = \@columndepth
ÊÊ \@linestogo = #1
ÊÊ \currentcolumn = 1
ÊÊ \def\@endcolumnactions{%
ÊÊÊÊÊÊ\ifnum \@linestogo<2
ÊÊÊÊÊÊÊÊ \the\crtok \egroup \endgroup \par 
ÊÊÊÊÊÊ\else
ÊÊÊÊÊÊÊÊ \global\advance\@linestogo by -1
ÊÊÊÊÊÊÊÊ \ifnum\@linestogoincolumn<2
ÊÊÊÊÊÊÊÊÊÊÊÊ\global\advance\currentcolumn by 1
ÊÊÊÊÊÊÊÊÊÊÊÊ\global\@linestogoincolumn = \@columndepth
ÊÊÊÊÊÊÊÊÊÊÊÊ\the\crtok
ÊÊÊÊÊÊÊÊ \else
ÊÊÊÊÊÊÊÊÊÊÊÊ&\global\advance\@linestogoincolumn by -1
ÊÊÊÊÊÊÊÊ \fi
ÊÊÊÊÊÊ\fi
ÊÊ }%
ÊÊ \makeactive\^^M
ÊÊ \letreturn \@endcolumnactions
ÊÊ \@columnwidth = \hsize
ÊÊÊÊ \advance\@columnwidth by -\parindent
ÊÊÊÊ \divide\@columnwidth by #2
ÊÊ \penalty\abovecolumnspenalty
ÊÊ \noindent 
ÊÊ \valign\bgroup
ÊÊÊÊ &\hbox to \@columnwidth{\strut \hsize = \@columnwidth ##\hfil}\cr
}%
\makeatother

\lefteqnumbers
   \def\testd{oui}
   \def\choixlat{\ifx\numadroite\testd\righteqnumbers
            \else  \lefteqnumbers\fi}
    \choixlat

\catcode`@=\letter
\def\@eplainfootnote#1{\let\@sf\empty 
  \ifhmode\edef\@sf{\spacefactor\the\spacefactor}\/\fi
  \global\advance\hlfootlabelnumber by 1
  \hlstart@impl{foot}{\hlfootlabel}%
  \hldest@impl{footback}{\hlfootbacklabel}%
  \hbox{$^{(#1)}$}%
  \hlend@impl{foot}%
  \@sf\vfootnote{#1.}%
}%
\catcode`@=\other

  \interfootnoteskip=0pt
  
  \everyfootnote={\eightpoint\leftskip=5truemm\rightskip5truemm}
  
  \hsize150truemm\vsize 240truemm\hoffset=5truemm
  
  \def\dimstand{\hsize 150truemm\vsize 250truemm\hoffset=8truemm\voffset=5truemm}
  
  \pretolerance=500\tolerance=1000\brokenpenalty=5000
  \parindent3mm
  
  \countdef\temps=170
  \temps=\time
  \countdef\nminutes=171{\nminutes = \time}
  \countdef\nheure=172
  \def\heure{\begingroup                     
     \temps = \time \divide\temps by 60
     \nheure = \temps                        
     \nminutes = \time
     \multiply\temps by 60
     \advance\nminutes by -\temps            
     \ifnum\nminutes<10 \toks1 = {0}%
     \else\toks1 = {}%
     \fi
     \number\nheure h\the\toks1 \number\nminutes  
  \endgroup}%

  \newcount\chstart
  \chstart=\pageno
 \headline={\ifnum\pageno=\chstart {\hfill} \else {\hss \tenrm --\ \folio\ --\hss}\fi}
  \footline={\hfill}
  \normalbaselines
  \frenchspacing
    \def\dater{\vglue-10mm\rightline{(\the\day/\the\month/\the\year)}}
  \def\dateheure{\vglue-10mm\rightline{(\the\day/\the\month/\the\year,\ \heure)}}

  \newif\ifpagetitre \pagetitretrue
\newtoks\hautpagetitre \hautpagetitre={\hfill}
\newtoks\baspagetitre \baspagetitre={\hfill}
\newtoks\auteurcourant \auteurcourant={\hfill}
\newtoks\titrecourant \titrecourant={\hfill}
\newtoks\hautpagegauche
\newtoks\hautpagedroite
\newtoks\hautpagemilieu
\hautpagemilieu={\tenrm\hfil -- \folio\ -- \hfil}
\hautpagegauche={\ifx\midfolio\oui\the\hautpagemilieu\else\tenrm\folio\hfill\the\auteurcourant\hfill\fi}
\hautpagedroite={\ifx\midfolio\oui\the\hautpagemilieu\else\hfill\the\titrecourant\hfill\tenrm\folio\fi}
\newtoks\baspagegauche \baspagegauche={\hfil}
\newtoks\baspagedroite \baspagedroite={\hfil}
\headline={\ifpagetitre\the\hautpagetitre
\else\ifodd\pageno\the\hautpagedroite\else\the\hautpagegauche\fi\fi }
\footline={\ifpagetitre\the\baspagetitre
\else\ifodd\pageno\the\baspagedroite
\else\the\baspagegauche\fi\fi \global\pagetitrefalse}

\def\pageblanche{\vfill\eject\pagetitretrue
\null\vfill\eject
\pagetitretrue
}
\def\chgtpage{\ifodd\pageno \else
\pageblanche \fi \pagetitretrue\titreun=0\footnotenumber=0}

\def\chgtpageincrtitreun{\ifodd\pageno \else
\pageblanche \fi \pagetitretrue\footnotenumber=0}

\def\majnombres{\ifodd\pageno \else
\pageblanche \fi \pagetitretrue\hautpoly\titreun=0\footnotenumber=0}

\def\hautspages#1#2{\auteurcourant={\ninepcap#1}\titrecourant={\nineit#2}}

\ifnum\chstart=\pageno \pagetitretrue\fi
  

  \def\PAR{\par}
  
  \def\leftnote#1{\vadjust{\setbox1=\vtop{\hsize 20mm\parindent=0pt\eightpoint
  \baselineskip=9pt\rightskip=4mm plus 4mm\vskip-4.7mm#1}\hbox{\kern-2cm\smash{\box1}}}}

  
  \def\raggedcenter{\leftskip=20pt plus 10em  
       \rightskip=\leftskip 
        \parfillskip=0pt 
         \spaceskip=.3333em \xspaceskip=.5em 
          \pretolerance=9999 \tolerance=9999
           \hyphenpenalty=9999 \exhyphenpenalty=9999 }
           
  \def\titrecentre#1{{\parindent0mm\raggedcenter
       \spaceskip=.6em plus .2em minus .2em\xspaceskip=.6em plus .2em minus .2em
        \tit#1\par}}
        


  \def\oui{oui}
  
   \def\fontetitreun{\ifx\paradouze\oui\douzepts\gpdouze\twelvebf\textfont1=\twelveib\else
\quatorzepts\gpquatorze\fourteenbf\fi}

\def\fontetitreunl{\douzepts\textfont1=\twelveib\scriptfont1=\tenib\fourteenti}
 
 \def\fontetitredeux{\textfont1=\eleveni\ifx\paradouze\oui\onzepts\scriptfont1=\ninei\elevenit\else
                        \douzepts\twelveit\fi}
 
   \def\fontetitredeuxb{\ifx\paradouze\oui\onzepts\eleventi\gponze\textfont1=\elevenib\scriptfont1=\nineib
                         \else\douzepts\twelveti\scriptfont1=\twelveib\scriptfont1=\tenib\gpdouze\fi}
                         
\def\fontetitredeuxl{\onzepts\textfont1=\elevenbf\scriptfont1=\ninebf\twelvebf}
  
\def\fontetitretrois{\textfont0=\elevenrm\scriptfont0=\eightrm\textfont1=\eleveni
                      \scriptfont1=\eighti\scriptscriptfont1=\sixi\elevenit}
                      
\def\fontetitrequatre{\textfont0=\elevenrm\scriptfont0=\eightrm\textfont1=\eleveni
                      \scriptfont1=\eighti\scriptscriptfont1=\sixi\elevenrm}
  
  \newcount\titreun\titreun=0
  \newcount\titredeux\titredeux=0
  \newcount\titretrois\titretrois=0
  \newcount\titrequatre\titrequatre=0
  \newcount\enonce\enonce=0
  
  \def\incr#1{\global\advance#1 by 1 {\the #1}}
  \def\avance#1{\global\advance#1 by 1}
  \def\init#1{\global#1=0}
  
  \long\def\Indentation#1#2{\setbox10=\hbox{\fontetitreun#1}
  	                    \ifdim\wd10 < 4mm
                         \setbox10=\hbox to 4mm{\box10\hfill}
                       \else\ifdim\wd10 < 6mm
                         \setbox10=\hbox to 6mm{\box10\hfill}
  	                    \else\ifdim\wd10 < 8mm
                         \setbox10=\hbox to 8mm{\box10\hfill}
                       \else\ifdim\wd10 < 12mm
                         \setbox10=\hbox to 12mm{\box10\hfill}
                       \fi\fi\fi\fi
                       \dimen10=\hsize
                       \advance \dimen10 by -\wd10
                       \noindent \box10 %
                       \ignorespaces
                       \hbox{\vtop{\hsize=\dimen10\raggedright\noindent\fontetitreun#2}}}

  \long\def\paraun#1{\removelastskip\par\medskip\goodbreak\vskip0pt plus.01\vsize\penalty-100
                \vskip0pt plus-.01\vsize
  	              \init{\titredeux}\ifnum\optionparag=1{\init\eqnumber\init\enonce}\else{}\fi
                  \goodbreak{\fontetitreun
  	                \Indentation{\incr{\titreun}.\ }{\fontetitreun #1\par}}\nobreak\medskip}

 %
 %
 \long\def\paraunc#1{\removelastskip\par\bigskip\goodbreak\vskip0pt plus.01\vsize\penalty-100
                \vskip0pt plus-.01\vsize
  	              \init{\titredeux}
                 \ifnum\optionparag=1{\init{\eqnumber}\init\enonce}\else{}\fi
                  \goodbreak
  	                {\parindent0mm\raggedcenter\fontetitreun\incr{\titreun}.\ 
                     \fontetitreun #1\par}\nobreak\medskip}
                     
\newtoks\titreunl
\titreunl={\ifnum\titreun=1{I}\fi%
\ifnum\titreun=2{II}\fi%
\ifnum\titreun=3{III}\fi%
\ifnum\titreun=4{IV}\fi%
\ifnum\titreun=5{V}\fi%
\ifnum\titreun=6{VI}\fi%
\ifnum\titreun=7{VII}\fi%
\ifnum\titreun=8{VIII}\fi%
\ifnum\titreun=9{IX}\fi%
\ifnum\titreun=10{X}\fi%
\ifnum\titreun=11{XI}\fi%
\ifnum\titreun=12{XII}\fi%
\ifnum\titreun=13{XIII}\fi%
}
\long\def\paraunl#1{\removelastskip\par\bigskip\bigskip\goodbreak\vskip0pt plus.01\vsize\penalty-100
                \vskip0pt plus-.01\vsize
  	              \init{\titredeux}\ifnum\optionparag=1{\init\eqnumber\init\enonce}\else{}\fi
                  \goodbreak{\fontetitreunl
  	                \Indentation{\global\advance\titreun by 1{\the\titreunl}.\ }{\fontetitreunl #1\par}}\nobreak\smallskip}

  
  \long\def\paradeux#1{\init{\titretrois}\vskip0pt plus.01\vsize\penalty-10
                \vskip0pt plus-.01\vsize\ifx \elie\oui\medskip\ifnum\titredeux>0\medskip\fi\fi
                 \Indentation{\fontetitredeux\the\titreun${\cdot}$\incr{\titredeux}.}
                              {\fontetitredeux\textfont1=\eleveni#1}\nobreak\par }
  
  \long\def\paradeuxb#1{\init{\titretrois}\vskip0pt plus.001\vsize\penalty-10
                \vskip0pt plus-.01\vsize{\ifx \elie\oui\medskip\ifnum\titredeux>0\medskip\fi\fi
                  \Indentation
  {\fontetitredeuxb\the\titreun${\cdot}$\incr{\titredeux}.}{ \fontetitredeuxb#1}}\nobreak
\smallskip}

\newtoks\titredeuxl
\titredeuxl={\ifnum\titredeux=1{A}\fi%
\ifnum\titredeux=2{B}\fi%
\ifnum\titredeux=3{C}\fi%
\ifnum\titredeux=4{D}\fi%
\ifnum\titredeux=5{E}\fi%
\ifnum\titredeux=6{F}\fi%
\ifnum\titredeux=7{G}\fi%
\ifnum\titredeux=8{H}\fi%
\ifnum\titredeux=9{I}\fi%
\ifnum\titredeux=10{J}\fi%
\ifnum\titredeux=11{K}\fi%
\ifnum\titredeux=12{L}\fi%
\ifnum\titredeux=13{M}\fi%
}
 \long\def\paradeuxl#1{\init{\titretrois}\vskip0pt plus.001\vsize\penalty-10
                \vskip0pt plus-.01
                \vsize \bigskip%
                  \Indentation
     {\fontetitredeuxl\global\advance\titredeux by 1
  \quad \the\titreunl${\cdot}$\the\titredeuxl.}{ \fontetitredeuxl#1}
  \removelastskip\nobreak\smallskip}
  

  \long\def\paratrois#1{\init{\titrequatre}\ifdim\lastskip<\smallskipamount
                \removelastskip\smallskip\fi
                 \vskip0pt plus.01\vsize\penalty-10
                  \vskip0pt
plus-.01\vsize{\ifx \elie\oui\ifnum\titretrois>0\medskip\fi\fi
\Indentation{\fontetitretrois\the\titreun${\cdot}$\the\titredeux${\cdot}$\incr{\titretrois}.\ }
  {\hskip0mm\baselineskip=14pt\fontetitretrois#1}\nobreak\smallskip}}
  
  
  \long\def\paratroisl#1{\init{\titrequatre}\ifdim\lastskip<\smallskipamount
                \removelastskip\fi
                 \vskip0pt plus.01\vsize\penalty-10
                  \vskip0pt
plus-.01\vsize\ifx \elie\oui\bigskip
\fi
\Indentation{\fontetitretrois\quad \quad \the\titreunl{${\cdot}$}\the\titredeuxl${\cdot}$\incr{\titretrois}.\ }
  {\hskip0mm\fontetitretrois#1}\nobreak\smallskip}


  \long\def\paraquatre#1{\ifdim\lastskip<\smallskipamount
                \removelastskip\smallskip\fi
                 \vskip0pt plus.01\vsize\penalty-10
                  \vskip0pt
                  plus-.01\vsize\par
 
\Indentation{\fontetitrequatre \the\titreun{${\cdot}$}\the\titredeux${\cdot}$\the\titretrois${\cdot}$\incr{\titrequatre}.\ }
{\hskip0mm\fontetitrequatre#1}\nobreak\smallskip}


\newtoks\titrequatrel
\titrequatrel={\ifnum\titrequatre=1{a}\fi%
\ifnum\titrequatre=2{b}\fi%
\ifnum\titrequatre=3{c}\fi%
\ifnum\titrequatre=4{d}\fi%
\ifnum\titrequatre=5{e}\fi%
\ifnum\titrequatre=6{f}\fi%
\ifnum\titrequatre=7{g}\fi%
\ifnum\titrequatre=8{h}\fi%
\ifnum\titrequatre=9{i}\fi%
\ifnum\titrequatre=10{j}\fi%
\ifnum\titrequatre=11{k}\fi%
\ifnum\titrequatre=12{l}\fi%
\ifnum\titrequatre=13{m}\fi%
}
\long\def\paraquatrel#1{\ifdim\lastskip<\smallskipamount
                \removelastskip\smallskip\fi
                 \vskip0pt plus.01\vsize\penalty-10
                  \vskip0pt
                  plus-.01\vsize{\bigskip
\Indentation{\global\advance\titrequatre by 1
\fontetitrequatre\quad \quad \quad \the\titreunl${\cdot}$\the\titredeuxl${\cdot}$\the\titretrois${\cdot}$\the\titrequatrel.\ }
{\hskip0mm\fontetitrequatre#1}\nobreak\smallskip}}

\ifx\optionkeys\oui
\def\drefun#1{\definexref{¤#1}{{\the\titreun}}{}} 
\def\drefdeux#1{\definexref{¤#1}{{\the\titreun}.{\the\titredeux}}{}}
\def\dreftrois#1{\definexref{¤#1}{{\the\titreun}.{\the\titredeux}.{\the\titretrois}}{}}
\else
\def\drefun#1{\definexref{prg#1}{{\the\titreun}}{}} 
\def\drefdeux#1{\definexref{prg#1}{{\the\titreun}.{\the\titredeux}}{}}
\def\dreftrois#1{\definexref{prg#1}{{\the\titreun}.{\the\titredeux}.{\the\titretrois}}{}}
\fi

%


  \long\def\propdeux#1#2#3#4{%
       \avance{\enonce}
       \leavevmode\edef\temp{#2}%
         \ifx\temp\empty 
          \else
           \definexref{#2}{#1~{\the\titreun.\the\enonce}}{enonces}
            \definexref{s#2}{{\the\titreun.\the\enonce}}{enonces}
             \fi
\smallskip
      \noindent{\bf#1\ {\bf\the\titreun.\the\enonce{#3}.}\enspace}{\sl#4\par}%
      \ifdim\lastskip<\medskipamount \removelastskip\penalty55\par \fi
   }

  \long\def\propun#1#2#3#4{%
      \avance{\enonce}
       \leavevmode\edef\temp{#2}%
        \ifx\temp\empty 
          \else
           \definexref{#2}{#1~{\the\enonce}}{enonces}
            \definexref{{s#2}}{{\the\enonce}}{enonces}
             \fi
   \par 
     \noindent{\bf#1\ {\bf\the\enonce{#3}.}\enspace}{\sl#4\par}%
     \ifdim\lastskip<\medskipamount \removelastskip\penalty55\medskip\fi
  }
  
  \long\def\prop#1#2#3#4{\ifnum\optionparag=1
                          \propdeux{#1}{#2}{\textfont1=\elevenib#3}{#4} \else\propun{#1}{#2}{\textfont1=\elevenib#3}{#4}\fi}

  \long\def\propt#1#2#3{\ifx\tpf\oui \prop{Th\'eo\-r\`eme}{#1}{#2}{#3}\par
                       \else\prop{Theorem}{#1}{#2}{#3}\par\fi}
  \long\def\Propt#1#2{\propt{#1}{}{#2}}
  \long\def\propl#1#2#3{\ifx\tpf\oui\prop{Lem\-me}{#1}{#2}{#3}\par
                         \else\prop{Lemma}{#1}{#2}{#3}\par\fi}
  
  \long\def\propc#1#2#3{\ifx\tpf\oui\prop{Corol\-laire}{#1}{#2}{#3}\par
                         \else\prop{Corollary}{#1}{#2}{#3}\par\fi}

  \long\def\propd#1#2#3{\ifx\tpf\oui\prop{D\'efi\-nition}{#1}{#2}{#3}\par
                       \else\prop{Definition}{#1}{#2}{#3}\par\fi} 
  
  \long\def\proptd#1#2#3{\ifx\tpf\oui\prop{Th\'eor\`eme et d\'efi\-nition}{#1}{#2}{#3}\par
                       \else\prop{Theorem and definition}{#1}{#2}{#3}\par\fi}


  
  \newcount\optionparag\optionparag=1
  
  \long\def\section#1#2{\ifnum\optionparag=1 \paraun{#2} 
                        \else\goodbreak{\fontetitreun
  	                \Indentation{#1.\ }{#2}}\nobreak\smallskip\fi}

  \def\eqconstruct#1{\ifnum\optionparag=1{\the\titreun\hbox{$\cdot$}#1}\else{#1}\fi}

  
  
  \def\numref{oui}  
  
  \newcount\mesref\mesref=0 
  \def\defbib#1{\ifx\numref\oui\global\advance\mesref by 1 \definexref{#1}{{\the
                 \mesref}}{}\else\definexref{#1}{#1}{}\fi}
  \def\bibtem#1{\defbib{#1}\item{\citer{#1}}}
  \def\citer#1{[\ref{#1}]}
  \def\citeplus#1#2{[\ref{#1}; #2]}

  
  \font\seventeenmsa=msam10 at 17pt    
  \font\fourteenmsa=msam10 at 14pt
  \font\twelvemsa=msam10 at 12pt
  \font\tenmsa=msam10                 
  \font\ninemsa=msam10 at 9pt 
  \font\eightmsa=msam10 at 8pt 
  \font\sevenmsa=msam7 
  \font\sixmsa=msam10 at 6pt
  \font\fivemsa=msam5
  \newfam\msafam\textfont\msafam=\tenmsa\scriptfont\msafam=\sevenmsa\scriptscriptfont\msafam=\fivemsa
  
  \font\seventeenbb=msbm10 at 17pt     
  \font\fourteenbb=msbm10 at 14pt
  \font\twelvebb=msbm10 at 12pt
  \font\tenbb=msbm10                   
  \font\ninebb=msbm10 at 9pt 
  \font\eightbb=msbm10 at 8pt 
  \font\sevenbb=msbm7 
  \font\sixbb=msbm10 at 6pt
  \font\fivebb=msbm5 
  \newfam\bbfam\textfont\bbfam=\tenbb\scriptfont\bbfam=\sevenbb\scriptscriptfont\bbfam=\fivebb
  \def\bb{\fam\bbfam\tenbb}%

  \font\seventeenscaln=eusm10 at 17pt   
  \font\twelvescaln=eusm10 at 12pt
  \font\tenscaln=eusm10                
  \font\ninescaln=eusm10 scaled 900
  \font\eightscaln=eusm10 scaled 800
  \font\sevenscaln=eusm10 scaled 700
  \font\sixscaln=eusm10 scaled 600
   
  \newfam\scalnfam\textfont\scalnfam=\tenscaln\scriptfont\scalnfam=\sevenscaln\scriptscriptfont\scalnfam=\sixscaln
  \def\scaln{\fam\scalnfam\tenscaln}%
  \def\scal{\scaln}
  
  \font\tenscalb=eusb10                

  \font\sevenscalb=eusb10 scaled 700

  \newfam\scalbfam\textfont\scalbfam=\tenscalb\scriptfont\scalbfam=\sevenscalb
  %
  
  %
  %
  \font\fourteenrm=cmr12 scaled 1200
  \font\elevenrm=cmr10 at 11pt
  \font\twelverm=cmr12
  \font\ninerm=cmr9
  \font\eightrm=cmr8      
  \font\sevenrm=cmr7
  \font\sixrm=cmr6

  \font\seventeenpcap=cmcsc10 at 17pt
  \font\tenpcap=cmcsc10                        
  \font\ninepcap=cmcsc9
  \font\eightpcap=cmcsc8
  \font\sevenpcap=cmcsc10 scaled 700
  
  \newfam\pcapfam\textfont\pcapfam=\tenpcap\scriptfont\pcapfam=\sevenpcap
  \def\pcap{\fam\pcapfam\tenpcap}
  
  \font\seventeenrm=cmbx12 scaled 1400

  \font\fourteenbf=cmbx10 scaled 1400
  
  \font\twelvebf=cmbx12
  \font\elevenbf=cmbx10 at 11pt
  \font\ninebf=cmbx9  
  \font\eightbf=cmbx8
  \font\sixbf=cmbx6
  
  \font\tengot=eufm10                           
   
  \font\eightgot=eufm10 at 8truept 
  \font\sevengot=eufm7 
  \font\sixgot=eufm10 at 6 truept 
   
  \newfam\gotfam
  \textfont\gotfam=\tengot\scriptfont\gotfam=\sevengot\scriptscriptfont\gotfam=\sixgot
  %

  
  \def\tit{%
  \textfont0=\seventeenrm\scriptfont0=\tenrm\def\rm{\fam0\seventeenrm}%
  \textfont1=\seventeenib\scriptfont1=\twelveib%
  \textfont2=\seventeensy\scriptfont2=\twelvesy\scriptscriptfont2=\ninesy
  \textfont3=\seventeenex
  \textfont\itfam=\seventeenti
  \def\it{\fam\itfam\seventeenti}%
  \textfont\bbfam=\seventeenbb \scriptfont\bbfam=\twelvebb
  \def\bb{\fam\bbfam\seventeenbb}%
  \textfont\msafam=\seventeenmsa\scriptfont\msafam=\twelvemsa
  \textfont\scalnfam=\seventeenscaln
  \def\pcap{\fam\pcapfam\seventeenpcap}
  \normalbaselineskip=25pt\normalbaselines\rm}

  \font\seventeenti=cmbxti10 scaled 1680
  
  \font\fourteenti=cmbxti10 at 14pt
  
  \font\twelveti=cmbxti10 scaled 1200
  \font\eleventi=cmbxti10 at 11pt

  %
  %
  \font\twelveit=cmti12	
  \font\elevenit=cmti10 scaled 1100
  \font\nineit=cmti9
  \font\eightit=cmti8
  \font\sevenit=cmti7

  %
  %
 
 \font\seventeenib=cmmib10 scaled 1680
  \font\fourteenib=cmmib10 scaled 1400
  \font\twelveib=cmmib10 scaled 1200
  \font\elevenib=cmmib10 scaled 1100
  \font\tenib=cmmib10
\font\eightib=cmmib10 scaled 800
  \font\nineib=cmmib10 scaled 900
\font\sevenib=cmmib10 scaled 700
\font\sixib=cmmib10 scaled 600
\font\fiveib=cmmib10 scaled 500

\ifx\ITAN\oui
\else
\innernewfam\cmmibfam
\textfont\cmmibfam=\tenib
\scriptfont\cmmibfam=\sevenib
\scriptscriptfont\cmmibfam=\fiveib
\def\ib{\fam\cmmibfam\tenib}
\fi

  %
  %
  \font\twelvei=cmmi10 scaled 1200
  \font\eleveni=cmmi10 scaled 1100
  \font\ninei=cmmi9
  \font\eighti=cmmi8 
  \font\seveni=cmmi7 			                
  \font\sixi=cmmi6
  
  \font\ninesl=cmsl9                    
  \font\eightsl=cmsl8 
  \font\sevensl=cmsl10 at 7pt

  \font\ninett=cmtt9                    
  \font\eighttt=cmtt8
  \font\seventt=cmtt10 scaled 700

  \font\seventeensy=cmsy10 scaled 1680    
  \font\fourteensy=cmsy10 scaled 1400
  \font\twelvesy=cmsy10 scaled 1176
  
  \font\ninesy=cmsy9                      
  \font\eightsy=cmsy8
  \font\sixsy=cmsy6
  \font\seventeenex=cmex10 at 17pt
  \font\fourteenex=cmex10 at 14pt
  \font\twelveex=cmex10 at 12pt
  \font\nineex=cmex10 at 9pt
  \font\eightex=cmex10 at 8pt
  \font\sevenex=cmex10 at 7pt
  \font\sixex=cmex10 at 6pt
  \font\fiveex=cmex10 at 5pt
  
   
  \font\fourteengp=cmmi10 at 14pt
  
  \font\twelvegp=cmmib10 at 12pt
  \font\elevengp=cmmib10 at 11pt
  \font\tengp=cmmib10                          
  \font\ninegp=cmmib10 at 9pt 
  \font\eightgp=cmmib8 
   
  \font\sixgp=cmmib6


  \def\gponze{\textfont0=\elevenbf\scriptfont0=\eightbf\scriptscriptfont0=\sixbf
           \textfont1=\elevengp\scriptfont1=\eightgp\scriptscriptfont1=\sixgp}
  \def\gpdouze{\textfont0=\twelvebf\scriptfont0=\tenbf\scriptscriptfont0=\ninebf
           \textfont1=\twelvegp\scriptfont1=\tengp\scriptscriptfont1=\ninegp}        
  
 \def\gpquatorze{\textfont0=\fourteenbf\scriptfont0=\twelvebf\scriptscriptfont0=\elevenbf
           \textfont1=\fourteengp\scriptfont1=\twelvegp\scriptscriptfont1=\elevengp}

  
  \expandafter\chardef\csname pre amssym.def at\endcsname=\the\catcode`\@
  \catcode`\@=11
  \def\undefine#1{\let#1\undefined}
  \def\newsymbol#1#2#3#4#5{\let\next@\relax
   \ifnum#2=\@ne\let\next@\msafam@\else
   \ifnum#2=\tw@\let\next@\bbfam@\fi\fi
   \mathchardef#1="#3\next@#4#5}
  \def\mathhexbox@#1#2#3{\relax
   \ifmmode\mathpalette{}{\m@th\mathchar"#1#2#3}%
   \else\leavevmode\hbox{$\m@th\mathchar"#1#2#3$}\fi}
  \def\hexnumber@#1{\ifcase#1 0\or 1\or 2\or 3\or 4\or 5\or 6\or 7\or 8\or
   9\or A\or B\or C\or D\or E\or F\fi}
  
  \def\setboxz@h{\setbox\z@\hbox}
  \def\wdz@{\wd\z@}
  \def\boxz@{\box\z@}
  
  \edef\msafam@{\hexnumber@\msafam}
  \mathchardef\dabar@"0\msafam@39
  
  \edef\bbfam@{\hexnumber@\bbfam}
  \def\widehat#1{\setboxz@h{$\m@th#1$}%
   \ifdim\wdz@>\tw@ em\mathaccent"0\bbfam@5B{#1}%
   \else\mathaccent"0362{#1}\fi}
  \def\widetilde#1{\setboxz@h{$\m@th#1$}%
   \ifdim\wdz@>\tw@ em\mathaccent"0\bbfam@5D{#1}%
   \else\mathaccent"0365{#1}\fi}
  \newsymbol\leqq 1335          
  \newsymbol\leqslant 1336
  \newsymbol\lessgtr 1337       
  \newsymbol\backprime 1038     
  \newsymbol\risingdotseq 133A  
  \newsymbol\fallingdotseq 133B 
  \newsymbol\succcurlyeq 133C   
  \newsymbol\geqq 133D          
  \newsymbol\geqslant 133E
  \newsymbol\nmid 232D
  \newsymbol\nexists 2040
  \newsymbol\smallsetminus 2272
  \newsymbol\varnothing 203F 
  \catcode`\@=\active

  \catcode`\@=11
  \newcount\typofr\typofr=1
  
  \catcode`\;=\active
  \def;{\ifnum\typofr=1\relax\ifhmode\ifdim\lastskip>\z@\unskip\fi
     \kern.2em\fi\string;\else\string;\fi}
  
  \catcode`\:=\active
  \def:{\ifnum\typofr=1\relax\ifhmode\ifdim\lastskip>\z@\unskip\fi
  \penalty\@M\ \fi\string:\else\string:\fi}
  
  \catcode`\!=\active
  \def!{\ifnum\typofr=1\relax\ifhmode\ifdim\lastskip>\z@\unskip\fi
     \kern.2em\fi\string!\else\string!\fi}
  
  \catcode`\?=\active
  \def?{\ifnum\typofr=1\relax\ifhmode\ifdim\lastskip>\z@\unskip\fi
     \kern.2em\fi\string?\else\string?\fi}

  \def\francais{\typofr=1\def\tpf{oui}}
  \def\anglais{\typofr=2\def\tpf{non}\def\english{oui}}
  \def\oui{oui}
  \francais
  
  \catcode`\@=12
  

\ifx\textures\oui
\def\raye #1|{\leavevmode\setbox1=\hbox{#1}%
\raise .5pt\hbox to \wd1{\xleaders\hbox{\rge{$ \sct / $}%
\kern 1pt}\hfill\kern -1pt }\kern -\wd1 \unhbox1\relax }

\def\barre #1|{\leavevmode\setbox1=\hbox{#1}%
\rlap{\Red\vrule height 2.4pt depth -1.2pt width \wd1}\Black \unhbox1\relax}
\else
\def\raye #1|{\leavevmode\setbox1=\hbox{#1}%
\raise .5pt\hbox to \wd1{\xleaders\hbox{\rge{$ \sct / $}%
\kern 1pt}\hfill\kern -1pt }\kern -\wd1 \unhbox1\relax }

\def\barre #1|{\leavevmode\setbox1=\hbox{#1}%
\rlap{\color{red}\vrule height 2.4pt depth -1.2pt width \wd1}\color{black} \unhbox1\relax}

\fi
  

  
  \def\og{\leavevmode\raise.24ex\hbox{$\scriptscriptstyle\langle\!\langle\>$}}    
  \def\fg{\leavevmode\raise.24ex\hbox{$\scriptscriptstyle\>\rangle\!\rangle$}}    

  \def\d{\,{\rm d}}

  \def\r{{\bb R}}
  
  \def\N{{\bb N}}

  \def\L{{\scal L}}
  \def\M{{\scal M}}
  
  \def\O{{\scal O}}
  \def\P{{\scaln P}}

  \def\frac#1#2{{#1\over #2}}
  \def\di#1#2{\sct#1\atop{\sct#2}}

  \def\numero{n$^{\rm o}\thinspace$}


  \def\numero{n$^{\rm o}\thinspace$}

  \def\¤{\S\thinspace}

  \def\¥{$\bullet$ }
  
  
  \def\e{{\rm e}}

  \def\epsilon{\varepsilon}

  \def\phi{\varphi}
  \def\theta{\vartheta}
  \def\rho{\varrho}
  \def\dm{{\textstyle{1\over 2}}}

  \def\sct{\scriptstyle}
  \def\pf{\noi{\it Proof. }}
  \def\nid{\ifnum\typofr=1\par\noindent{\it D\'emonstration. }\else\pf\fi}
  \def\noi{\noindent}
  \def\rem{\ifnum\typofr=1\noi{\it Remarque.}\ \else\noi{\it Remark.}\ \fi}
  \def\rems{\ifnum\typofr=1\noi{\it Remarques.}\ \else\noi{\it Remarks.}\ \fi}
  \def\re{{\Re e\,}}

  \def\1{{\bf 1}}
  \def\|{\Vert}

  \def\leq{\leqslant}
  \def\geq{\geqslant}

  \def\eg{{e.g.}}
  \newsymbol\subsetneqq 2324
  \newsymbol\subsetneq 2328

  \def\fl#1{\left\lfloor #1 \right\rfloor}

  \def\log{\mathop{\rm log}\nolimits}
  
  \def\fs#1#2{{\scriptstyle{#1\over #2}}}
  


\def\abs#1{\left|#1\right|}


  \def\pmb#1{\setbox0=\hbox{#1}%
  \kern-.025em\copy0\kern-\wd0\kern.05em\copy0\kern-\wd0\kern-.025em\raise .0433em\box0 }

  
  \skewchar\eighti='177 \skewchar\sixi='177
  \skewchar\eightsy='60 \skewchar\sixsy='60
  
  \def\eightpoint{%
  \textfont0=\eightrm\scriptfont0=\sixrm\scriptscriptfont0=\fiverm
  \def\rm{\fam0\eightrm}%
  \textfont1=\eighti\scriptfont1=\sixi
  \scriptscriptfont1=\fivei\def\oldstyle{\fam1\seveni}%
  \textfont2=\eightsy\scriptfont2=\sixsy\scriptscriptfont2=\fivesy
  \textfont3=\eightex\scriptfont3=\sixex
  \textfont\itfam=\eightit
  \def\it{\fam\itfam\eightit}%
  \textfont\slfam=\eightsl
  \def\sl{\fam\slfam\eightsl}%
  \textfont\bbfam=\eightbb \scriptfont\bbfam=\sixbb\scriptscriptfont\bbfam=\fivebb
  \def\bb{\fam\bbfam\eightbb}%
  \textfont\msafam=\eightmsa\scriptfont\msafam=\sixmsa
  \textfont\scalnfam=\eightscaln
  \def\scaln{\fam\scalnfam\eightscaln}
  \textfont\ttfam=\eighttt
  \def\tt{\fam\ttfam\eighttt}%
\textfont\gotfam=\eightgot
  \textfont\bffam=\eightbf\scriptfont\bffam=\sixbf\scriptscriptfont\bffam=\fivebf
  \def\bf{\fam\bffam\eightbf}%
  \ifx\ITAN\oui\else\textfont\cmmibfam=\eightib
       \scriptfont\cmmibfam=\sixib
        \scriptscriptfont\cmmibfam=\fiveib
         \def\ib{\fam\cmmibfam\eightib}
   \fi
  \textfont\pcapfam=\eightpcap
  \def\pcap{\fam\pcapfam\eightpcap}
  \abovedisplayskip=2pt plus2pt minus 2pt
  \belowdisplayskip=2pt plus1pt minus 2pt
  \abovedisplayshortskip= 1pt plus 2pt minus 1pt
  \belowdisplayshortskip= 1pt plus 2pt minus 1pt
  \smallskipamount=2pt plus 1pt minus 2pt
  \medskipamount=3pt plus 2pt minus 2pt
  \bigskipamount=7pt plus 3pt minus 3pt
  \setbox\strutbox=\hbox{\vrule height 5pt depth 2pt width 0pt}%
  \normalbaselineskip=9pt\normalbaselines\rm}

  \def\({\left(}
  \def\){\right)}
  
  \def\footnoterule{\kern -2pt\hrule width 7truecm\kern 2.4pt}
  
  \def\xnotedef#1{\definexref{#1}{\noexpand\number\footnotenumber}{Note}}%

  
  
  \def\ninepoint{%
  \textfont0=\ninerm\scriptfont0=\sixrm\scriptscriptfont0=\fiverm
  \def\rm{\fam0\ninerm}%
  \textfont1=\ninei\scriptfont1=\sixi
  \scriptscriptfont1=\fivei\def\oldstyle{\fam1\ninei}%
  \textfont2=\ninesy\scriptfont2=\sixsy\scriptscriptfont2=\fivesy
  \textfont3=\nineex\scriptfont3=\sixex
  \textfont\itfam=\nineit
  \def\it{\fam\itfam\nineit}%
  \textfont\slfam=\ninesl
  \def\sl{\fam\slfam\ninesl}%
  \textfont\bbfam=\ninebb\scriptfont\bbfam=\sixbb\scriptscriptfont\bbfam=\fivebb
  \def\bb{\fam\bbfam\ninebb}%
  \textfont\msafam=\ninemsa\scriptfont\msafam=\sixmsa\scriptscriptfont\msafam=\fivemsa
  \textfont\scalnfam=\ninescaln
  \def\scaln{\fam\scalnfam\ninescaln}
  \textfont\ttfam=\ninett
  \def\tt{\fam\ttfam\ninett}%
  \textfont\bffam=\ninebf\scriptfont\bffam=\sixbf\scriptscriptfont\bffam=\fivebf
  \def\bf{\fam\bffam\ninebf}%
  \abovedisplayskip=3pt plus2pt minus 2pt
  \belowdisplayskip=3pt plus1pt minus 2pt
  \abovedisplayshortskip= 2pt plus 2pt minus 1pt
  \belowdisplayshortskip= 2pt plus 2pt minus 1pt
  \smallskipamount=2pt plus 1pt minus 2pt
  \medskipamount=3pt plus 2pt minus 2pt
  \bigskipamount=7pt plus 3pt minus 3pt
  \setbox\strutbox=\hbox{\vrule height 5pt depth 2pt width 0pt}%
  \normalbaselineskip=11pt plus.3pt minus.2pt\normalbaselines\rm}

  \def\sevenpoint{%
  \textfont0=\sevenrm\scriptfont0=\sixrm\scriptscriptfont0=\fiverm
  \def\rm{\fam0\sevenrm}%
  \textfont1=\seveni\scriptfont1=\sixi
  \scriptscriptfont1=\fivei\def\oldstyle{\fam1\seveni}%
  \textfont2=\sevensy\scriptfont2=\sixsy\scriptscriptfont2=\fivesy
  \textfont3=\sevenex\scriptfont3=\fiveex
  \textfont\itfam=\sevenit
  \def\it{\fam\itfam\sevenit}%
  \textfont\slfam=\sevensl
  \def\sl{\fam\slfam\sevensl}%
  \textfont\bbfam=\sevenbb \scriptfont\bbfam=\sixbb\scriptscriptfont\bbfam=\fivebb
  \def\bb{\fam\bbfam\sevenbb}%
  \textfont\msafam=\sevenmsa\scriptfont\msafam=\sixmsa
  \textfont\scalnfam=\sevenscaln
  \def\scaln{\fam\scalnfam\sevenscaln}
  \textfont\bffam=\sevenbf\scriptfont\bffam=\sixbf\scriptscriptfont\bffam=\fivebf
  \def\bf{\fam\bffam\sevenbf}%
  \textfont\ttfam=\seventt
  \abovedisplayskip=2pt plus2pt minus 2pt
  \belowdisplayskip=2pt plus1pt minus 2pt
  \abovedisplayshortskip= 1pt plus 2pt minus 1pt
  \belowdisplayshortskip= 1pt plus 2pt minus 1pt
  \smallskipamount=2pt plus 1pt minus 2pt
  \medskipamount=3pt plus 2pt minus 2pt
  \bigskipamount=7pt plus 3pt minus 3pt
  \setbox\strutbox=\hbox{\vrule height 5pt depth 2pt width 0pt}%
  \normalbaselineskip=9pt\normalbaselines\rm}

 \def\onzepts{%
 \textfont0=\elevenrm\scriptfont0=\ninerm
 \textfont1=\eleveni\scriptfont1=\ninei
}

\def\douzepts{%
  \textfont0=\twelverm\scriptfont0=\tenrm\def\rm{\fam0\twelverm}%
  \textfont1=\twelvei\scriptfont1=\teni%
  \textfont2=\twelvesy\scriptfont2=\tensy\scriptscriptfont2=\eightsy
  \textfont3=\twelveex
  \textfont\itfam=\twelveti
  \def\it{\fam\itfam\twelveti}%
  \textfont\bffam=\twelvebf\scriptfont\bffam=\tenbf\scriptscriptfont\bffam=\eightbf
  \def\bf{\fam\bffam\twelvebf}%
  \textfont\bbfam=\twelvebb \scriptfont\bbfam=\tenbb
  \def\bb{\fam\bbfam\twelvebb}%
  \textfont\msafam=\twelvemsa\scriptfont\msafam=\tenmsa
  \textfont\scalnfam=\twelvescaln
  \normalbaselineskip=15pt\normalbaselines\rm}

\def\quatorzepts{%
  \textfont0=\fourteenrm\scriptfont0=\twelverm\def\rm{\fam0\fourteenrm}%
  \textfont1=\fourteenib\scriptfont1=\twelveib%
  \textfont2=\fourteensy\scriptfont2=\twelvesy\scriptscriptfont2=\tensy
  \textfont3=\fourteenex
  \textfont\itfam=\fourteenti
  \def\it{\fam\itfam\fourteenti}%
  \textfont\bffam=\fourteenbf\scriptfont\bffam=\twelvebf\scriptscriptfont\bffam=\tenbf
  \def\bf{\fam\bffam\fourteenbf}%
  \textfont\bbfam=\fourteenbb \scriptfont\bbfam=\twelvebb
  \def\bb{\fam\bbfam\fourteenbb}%
  \textfont\msafam=\fourteenmsa\scriptfont\msafam=\twelvemsa
  \textfont\scalnfam=\twelvescaln
  \normalbaselineskip=18pt\normalbaselines\rm}


\def\AA{{\it Acta Arith.}}

\def\picture #1 by #2 (#3){\leavevmode\vbox to #2{
     \hrule width #1 height 0pt depth 0pt
      \vfill
       \special{picture #3}}}

\def\illustration #1 by #2 (#3) scaled #4{\dimen1=#2
  \divide\dimen1 by 1000
  \multiply\dimen1 by #4
  \vtop to \dimen1{\dimen1=#1
  \divide\dimen1 by 1000
  \multiply\dimen1 by #4
  \hsize=\dimen1\vss
  \noindent\special{illustration #3 scaled #4}}}

\ifx\couleurs\oui

\fi

 \fi
\newcount\numchapi\numchapi=0

\ifx\optionkeymacros\undefined\else \fi

\catcode`\Œ=\active\defŒ{{\aa}}       
\catcode`\º=\active\defº{\int}        
\catcode`\=\active\def{\c c}        
\catcode`\¶=\active\def¶{\partial}    
\catcode`\Ä=\active\defÄ{\oint}       
\catcode`\Æ=\active\defÆ{\triangle}   
\catcode`\Â=\active\defÂ{\neg}        
\catcode`\µ=\active\defµ{\mu}         
\catcode`\¿=\active\def¿{{\o}}        
\catcode`\¹=\active\def¹{\pi}         
\catcode`\Ï=\active\defÏ{{\oe}}       
\catcode`\§=\active\def§{{\ss}}       
\catcode`\ =\active\def {\dagger}     
\catcode`\Ã=\active\defÃ{\sqrt}       
\catcode`\·=\active\def·{\Sigma}      
\catcode`\Å=\active\defÅ{\approx}     
\catcode`\½=\active\def½{\Omega}      
\catcode`\£=\active\def£{{\it\$}}     
\catcode`\°=\active\def°{\infty}      
\catcode`\¤=\active\def¤{{\S}}        
\catcode`\¦=\active\def¦{{\P}}        
\catcode`\¥=\active\def¥{\bullet}     
\catcode`\»=\active\def»{\leavevmode\raise.585ex\hbox{\b a}}      
\catcode`\¼=\active\def¼{\leavevmode\raise.6ex\hbox{\b o}}        
\catcode`\­=\active\def­{\not=}       
\catcode`\²=\active\def²{\leq}        
\catcode`\³=\active\def³{\geq}        
\catcode`\Ö=\active\defÖ{\div}        
\catcode`\É=\active\defÉ{{\dots}}     
\catcode`\¾=\active\def¾{{\ae}}       
\catcode`\Ç=\active\defÇ{\og}         
\catcode`\Ò=\active\defÒ{``}          
\catcode`\Á=\active\defÁ{!`}          
\catcode`\¢=\active\def¢{\rlap/c}     
\catcode`\Ô=\active\defÔ{`}           
\catcode`\Õ=\active\defÕ{'}           


\catcode`\=\active\def{{\AA}}       
\catcode`\'=\active\def'{\c C}        
\catcode`\¯=\active\def¯{{\O}}        
\catcode`\¸=\active\def¸{\Pi}         
\catcode`\Î=\active\defÎ{{\OE}}       
\catcode`\®=\active\def®{{\AE}}       
\catcode`\×=\active\def×{\diamond}    
\catcode`\¡=\active\def¡{\accent'27}  
\catcode`\Ó=\active\defÓ{''}          
\catcode`\±=\active\def±{\pm}         
\catcode`\È=\active\defÈ{\fg}         
\catcode`\À=\active\defÀ{?`}          
\catcode`\Ð=\active\defÐ{--}          
\catcode`\Ñ=\active\defÑ{---}         


\catcode`\Š=\active\defŠ{\"a}        
\catcode`\'=\active\def'{\"e}        
\catcode`\•=\active\def•{\"{\i}}     
\catcode`\š=\active\defš{\"o}        
\catcode`\Ÿ=\active\defŸ{\"u}        
\catcode`\Ø=\active\defØ{\"y}        
\catcode`\å=\active\defå{\^A}        
\catcode`\€=\active\def€{\"A}        
\catcode`\…=\active\def…{\"O}        
\catcode`\†=\active\def†{\"U}        
\catcode`\‡=\active\def‡{\'a}        
\catcode`\Ž=\active\defŽ{\'e}        
\catcode`\'=\active\def'{\'{\i}}     
\catcode`\—=\active\def—{\'o}        
\catcode`\œ=\active\defœ{\'u}        
\catcode`\ƒ=\active\defƒ{\'E}        
\catcode`\æ=\active\defæ{\^E}        
\catcode`\é=\active\defé{\`E}        %
\catcode`\ˆ=\active\defˆ{\`a}        
\catcode`\=\active\def{\`e}        
\catcode`\"=\active\def"{\`{\i}}     
\catcode`\˜=\active\def˜{\`o}        
\catcode`\=\active\def{\`u}        
\catcode`\Ë=\active\defË{\`A}        
\catcode`\‹=\active\def‹{\~a}        
\catcode`\–=\active\def–{\~n}        
\catcode`\›=\active\def›{\~o}        
\catcode`\Ì=\active\defÌ{\~A}        
\catcode`\"=\active\def"{\~N}        
\catcode`\Í=\active\defÍ{\~O}        
\catcode`\‰=\active\def‰{\^a}        
\catcode`\=\active\def{\^e}        
\catcode`\"=\active\def"{\^{\i}}     
\catcode`\™=\active\def™{\^o}        
\catcode`\ž=\active\defž{\^u}        

\let\optionkeymacros\null

\optionparag=2
\def\paradouze{oui}
\anglais
\def\FF{{\scal F}}
\dimstand
\voffset-3mm
\vsize=255truemm
\hautspages{R. de la Bretche \& G. Tenenbaum}{On the gap distribution of prime factors}
\dateheure

\titrecentre{On the gap distribution of prime factors}
\bigskip\medskip
\centerline{R. de la Bretche \& G. Tenenbaum}
\bigskip\bigskip
{\eightpoint\leftskip1cm\rightskip1cm
\noi{\bf Abstract.} Let $\{p_j(n)\}_{j=1}^{\omega(n)}$ denote the increasing sequence of distinct prime factors of an integer~$n$. 
For $z\geqslant 0$, let  $G(n;z)$ denote the number of those indexes $j$ such that $p_{j+1}(n)>p_j(n)^{\exp z}$. We show uniform convergence, with almost optimal effective estimate of the speed, of the distribution of $G(n;z)$ on $\{n:1\leqslant n\leqslant N\}$ to a Gaussian limit law with mean $\e^{-z}\log_2n$ and variance $\{\e^{-z}-2z\e^{-2z}\}\log_2n$, and we establish an asymptotic formula with remainder for all centered moments.  \PAR
\medskip\noi
{ \bf Keywords:} Distribution of prime factors, normal order, gaps between prime factors, moments of arithmetical functions, Fourier inversion, Laplace transform, Berry-Esseen inequality. \par
\smallskip 
\noi \bf 2020 Mathematics Subject Classification: \rm  11N25, 11N35.\par }
\par 
\bigskip\bigskip
\paraun{Introduction}
Let $\{p_j(n)\}_{j=1}^{\omega(n)}$ denote the increasing sequence of distinct prime  factors of an integer $n$. For $z>0$, define
$$
G(n;z):=|\{1\leqslant j<\omega(n):p_{j+1}(n)>p_j(n)^{\exp z}\}|\qquad (n\geqslant 1,\,z\geqslant 0),\quad\lambda=\lambda_z:=\e^{-z}.$$
By the sieve, it is easy to see that, for  fixed $z>0$, $N\to\infty$, and $j=j_N\to\infty$, we have $p_{j+1}(n)>p_j(n)^{\exp z}$ for $\lambda N+o(N)$ integers $n\leqslant N$, and so that the average order of $G(n;z)$ is $\{\lambda+o(1)\}\log_2n$. Here and in the sequel, we denote by $\log_k$ the $k$th fold iterated logarithm, defined for sufficient large value of the argument.
\par 
Write $G^*_N(n;z):=G(n,z)-\lambda\log_22N$ $(1\leqslant n\leqslant N)$. Erd\H os~\citer{Er46} showed that $G^*_N(n;z)^2$ has mean-value $o\big((\log_2N)^2\big)$ over the set of integers $n\leqslant N$,  thus deducing that $\lambda\log_2n$ is a normal order for $G(n;z)$. In the same spirit,  he  stated later \citer{Er59} that, for any given $\lambda>0$, the inequality $\max_{1\leqslant j<\omega(n)}(\log p_{j+1}(n)\}/\log p_j(n)>(\log_2n)^\lambda$ holds on a sequence of integers $n$ with asymptotic density~$1-\lambda$.  The second author \citer{Te19}  provided the necessary details for the proof of this statement.\par 
In a recent preprint, Sofos \citer{So21} studied the centered moments of $\{G(n;z)\}_{1\leqslant n\leqslant\leqslant N}$
and established, for each fixed $z\geqslant 0$, the asymptotic estimate
$$\sum_{n\leqslant N}G^*_N(n;z)^r=\{\mu_r+o(1)\}N(\gamma\log_2N)^{r/2}\qquad (r\in\N,\,N\to\infty)$$
with $\gamma=\gamma_z:=\lambda(1-2\lambda z)$, where $\mu_r$ denotes the $r$th moment of the normal law. 
\par 
Sofos' method rests on a probabilistic theorem of Stein \citer{St86}, which provides a quantitative estimate for the speed of weak convergence of a sum of weakly dependent random variables to the Gaussian law. A second, crucial part of the analysis   involves  studying the convergence of moments of the model. The proof is then completed by showing sufficient agreement between the model and the arithmetical setting.    According to the terminology introduced by Ruzsa \citer{Ru82}, this proof thus pertains to the class of {\it indirect} approaches, in which arithmetic results are derived from  genuine probabilistic statements via  comparison theorems between the set $\{n:1\leqslant n\leqslant N\}$ $(N\in\N)$ equipped with  uniform probability and a suitable probabilistic model. 
\par 
Here, we propose a {\it direct} approach in which the Fourier transform
$${1\over N}\sum_{n\leqslant N}\e^{iyG^*_N(n;z)}\qquad (y\in\r,\,N\geqslant 1)\eqdef{TL}$$
is evaluated by purely arithmetic means, with  sufficient accuracy to apply an effective  inversion result such as the Berry-Esseen inequality---see, \eg, \citeplus{Te15}{th. II.7.16}. An upper bound for the corresponding bilateral Laplace transform then readily solves, in a comparable effective way, the question of convergence of moments.
This method has three advantages: (i) the proof turns out to be very short;  (ii) the speed of convergence is  explicitly bounded; (iii) the effective asymptotic formula for  moments is a straight corollary of the Fourier inversion theorem and the Laplace transform upper bound. 
\par \medskip\goodbreak
Writing $\Phi(t):=(1/\sqrt{2\pi})\int_{-\infty}^t\e^{-y^2/2}\d y$ $(t\in \r)$, we prove the following result.\par 
\Propt{LoiG}{Let $z\geqslant 0$, $\lambda_z:=\e^{-z}$, $\gamma_z:=\lambda_z(1-2z\lambda_z)$. Uniformly for $t\in\r$, $N\geqslant 16$, we have
$${1\over N}\sum_{\di{n\leqslant N}{G^*_N(n;z)\leqslant t\sqrt{\gamma_z\log_2N}}}1=\Phi(t)+O\bigg({\log_3N\over \sqrt{\log_2 N}}\bigg).\eqdef{floiG}$$ 
Moreover, for any given non-negative integer $r$, we have
$${1\over N}\sum_{n\leqslant N}G^*_N(n;z)^{r}=(\gamma_z\log_2N)^{r/2}\bigg\{\mu_r+O\bigg({(\log_3N)^{1+r/2}\over \sqrt{\log_2 N}}\bigg)\bigg\},\eqdef{evalmoms}$$
where $\mu_{r}$ is the $r$th moment of the normal law. This formula persists for all real $r\geqslant 0$ provided $G^*_N(n;z)^r$ is replaced by $|G^*_N(n;z)|^r$ and $\mu_r$ by $2^{r/2}\Gamma(\dm r+\dm)/\sqrt{\pi}$.}
As follows from a sharp version \citer{De59} of the Erd\H os-Kac theorem, which corresponds to $z=0$,   the error terms in \eqref{floiG} and \eqref{evalmoms} are optimal up to the factor involving $\log_3N$. 
\paraun{Proof}
Assuming throughout that $z>0$, we first prove \eqref{floiG}. Let $\chi_p(n)$ denote the indicator  of the set of those integers $n\leqslant N$ 
that are divisible by the prime $p$ but by no prime $q$ such that $p< q\leqslant  p^{1/\lambda}$. 
\par 
Let $T=T_N$ to be precisely defined later, with $T_N=o\big((\log_2N)/\log_4N\big)$. Put $$u:=\exp\{(\log N)^{1/T}\},\quad v:=\exp\{(\log N)^{1-1/T}\},\qquad w:=\log \Big({\log v\over \log u}\Big),$$
and select $T$ in such a way that $w/z$ is an integer. We define
$$g(n):=\sum_{u<p\leqslant v}\chi_p(n)-\lambda w,$$
and aim to show that the distribution of $g(n)/\sqrt{w}$  is approximately Gaussian over the set of integers $n\in[1,N]$. With an effective Fourier inversion in mind, we consequently introduce the quantity
$$L(N;\vartheta):=\sum_{n\leqslant N}\e^{i\vartheta g(n)}\qquad (\vartheta\in\r),$$
with $\vartheta:=y/\sqrt{w}$. We shall only need to consider  $|y|\ll \sqrt{w}$, hence $|\vartheta|\ll 1$.\par  We have
$$L(N;\vartheta)=\e^{-i\vartheta \lambda w}\sum_{n\leqslant N}\prod_{p|n}\big\{1+\chi_p(n)(\e^{i\vartheta}-1)\big\}.$$
Let $\M=\M(u,v)$ denote the set of those squarefree integers $m\leqslant N$ all of whose prime factors belong to the interval $]u, v]$ and such that $\chi_p(m)=1$ whenever $p|m$. For each $m\in\M$, define $P_m:=\prod_{p|m}\prod_{p< q\leqslant p^{1/\lambda}}q$. Thus
$$L(N;\vartheta)=\e^{-i\vartheta \lambda w}\sum_{m\in\M}(\e^{i\vartheta}-1)^{\omega(m)}\sum_{\di{d\leqslant N/m}{(d,P_m)=1}}1.\eqdef{Lxt1}$$
A standard application of Rankin's method yields that the contribution of those $m>\sqrt{N}$ is $\ll N/u$.  By a sieve result, e.g. \citeplus{Te93}{lemma}, that of the remaining integers equals
$${N\varphi(P_m)\over mP_m}\Big\{1+O\Big({1\over u}\Big)\Big\}={N\lambda^{\omega(m)}\over m}\exp\Big\{O\Big({1+\omega(m)\over \exp\sqrt{\log u}}\Big)\Big\},$$
where the second estimate follows from a strong form of the prime number theorem.
Observing that, by definition of $\M$ we must have $z\omega(m)\leqslant w$, we obtain
$$L(N;\vartheta)=N\e^{-i\vartheta \lambda w}\sum_{k\leqslant w/z}\lambda^k(\e^{i\vartheta}-1)^k\sum_{\di{m\in\M}{\omega(m)=k}}{1\over m}+O_z\Big(N\e^{-\fs12\sqrt{\log u}}\Big).$$
 In the inner sum, we may write $m=p_1\cdots p_k$, with $$\log u<\log p_1\leqslant \lambda\log p_2\leqslant \lambda^2\log p_3\leqslant \cdots\leqslant \lambda^{k-1}\log p_k\leqslant \lambda^{k-1}\log v.$$
The sum over $p_k$ is
$$\log_2v-z-\log_2p_{k-1}+O\Big({1\over \log p_{k-1}}\Big).$$ 
Summing this over $p_{k-1}$ furnishes in turn
$$\dm(\log_2v-2z-\log_2p_{k-2})^2+O\Big({1\over \log p_{k-2}}\Big).$$
Iterating, the sum over $p_{k-j-1}$ is
$${1\over j!}(\log_2v-jz-\log_2p_{k-j})^j+O\Big({1\over \log p_{k-j}}+{j\over (\log p_{k-j})^2}\Big)$$
and finally
$$\sum_{\di{m\in\M}{\omega(m)=k}}{1\over m}={(w-kz)^k\over k!}+O\Big({1\over \log u}\Big).$$
Therefore
$$L(N;\vartheta)=NF(\vartheta,w)+O\Big({N\over (\log N)^{1/T}}\Big),\eqdef{Lxt2}$$
with $$F(\vartheta,w):=\e^{-i\lambda\vartheta w}\sum_{k\leqslant w/z}{\lambda^k(\e^{i\vartheta}-1)^k(w-kz)^k\over k!}.\eqdef{defF}$$
\par 
Write $\Delta(\vartheta,w):=F(\vartheta,w)-\e^{-\fs12\vartheta^2\gamma w}$. We now aim to show that, still with $\gamma=\gamma_z$ and for suitable positive constants $c_0$, $c_1$ depending at most on $z$,  we have
$$\leqalignno{
\Delta(\vartheta,w)&\ll \vartheta\qquad (|\vartheta|\leqslant w^{-2}),&\eqdef{eval1F}\cr
\Delta(\vartheta,w)&\ll \vartheta^2w\qquad (w^{-2}< |\vartheta|\leqslant w^{-3/4}),&\eqdef{eval2F}\cr
\Delta(\vartheta,w)&\ll \vartheta^3w\e^{-\fs12\gamma\vartheta^2w}\qquad (w^{-3/4}< |\vartheta|\leqslant w^{-1/3}),&\eqdef{evalF}\cr
\Delta(\vartheta,w)&\ll \e^{-c_1\vartheta^2w}\qquad (w^{-1/3}<|\vartheta|\leqslant c_0),&\eqdef{eval3F}}$$
\par 
We first note that \eqref{eval1F} readily follows from the estimates
$$\e^{-i\lambda\vartheta w}=1-i\lambda\vartheta w+O(\vartheta),\quad \sum_{k\leqslant w/z}{\lambda^k(\e^{i\vartheta}-1)^k(w-kz)^k\over k!}=1+i\lambda \vartheta w+O(\vartheta)\quad (|\vartheta|\leqslant w^{-2}).$$
\par 
We next prove \eqref{eval2F}. We introduce the notation
$$\mu:=\lambda(\e^{i\vartheta}-1), \quad \varrho:=|\mu|=2\lambda|\sin(\dm\vartheta)|.\eqdef{notmu}$$
Put $m:=\fl{\varrho w+\sqrt{w}}.$ Then
$$\sum_{k\geqslant m+1}{\varrho^kw^k\over k!}\ll {\varrho^mw^m\over m!}\ll\Big({\varrho\e w\over m}\Big)^m\ll\e^{-\sqrt{w}}\ll\vartheta^2w.$$
Therefore, using the fact that $(w-kz)^{k}=w^k\{1+O(k^2/w)\}$ for $k\leqslant m$,
$$F(\vartheta,w)=\e^{-i\lambda\vartheta w}\Big\{\e^{\mu w}+O(\vartheta^2w)\Big\}=\e^{(\mu-i\lambda \vartheta)w}+O(\vartheta^2 w)=1+O(\vartheta^2 w).$$
This is plainly sufficient.
\par 
We now embark to proving \eqref{evalF} and \eqref{eval3F}.
  Define $K:=w/z$ (an integer by construction), so that, for any $R>0$,
  $$\eqalign{\e^{i\lambda\vartheta w}F(\vartheta,w)&=\sum_{0\leqslant k\leqslant K}{\mu^k(w-kz)^k\over k!}={1\over 2\pi i}\oint_{|\zeta|=R}\sum_{0\leqslant k\leqslant K}{\e^{\mu (w-kz)\zeta}\over \zeta^{k+1}}\d\zeta\cr
  &={1\over 2\pi i}\oint_{|\zeta|=R}{1-\e^{-(K+1)\mu z\zeta}/\zeta^{K+1}\over \zeta-\e^{-\mu z\zeta}}\e^{\mu w\zeta}\d\zeta.\cr}$$
Selecting $R:=1/|\vartheta|$, we see that the contribution of the term $\e^{-(K+1)\mu z\zeta}/\zeta^{K+1}$ is $$\ll \e^{2w}|\vartheta|^{w/z}\ll \e^{-w}$$
in the considered range. In the remaining integral, we observe that $-\re(\mu z\zeta)\ll1$, so, by RouchŽ's theorem, the equation $\zeta-\e^{-\mu z\zeta}=0$ has exactly one solution, say $\beta$, inside the circle. Moreover, we have $\beta=1-\mu z+O(\vartheta^2)$.  The residue theorem thus yields
$$\e^{i\lambda\vartheta w}F(\vartheta,w)={\e^{\mu w\beta}\over 1+\mu z\e^{-\mu z\beta}}+O\big(\e^{-w}\big)=\e^{\mu w-\mu^2zw}\{1+O(\vartheta)\}\quad(w^{-3/4}<|\vartheta|\leqslant w^{-1/3}).$$
Since $-i\lambda\vartheta+\mu-\mu^2z=-\dm\gamma\vartheta^2+O(\vartheta^3)$, this completes the proof of \eqref{evalF}. When $x^{-1/3}<|\vartheta|\leqslant c_0$ with sufficiently small $c_0=c_0(z)$, we note that $\re\{\mu w\beta-i\lambda\vartheta w\}\leqslant -c_1\vartheta^2w$: this confirms~\eqref{eval3F}.
\par 
Selecting with $T_N:=  \big\{1+O(1/\log_3N)\big\}(\log_2N)/\log_3N$, we may now apply the Berry-Esseen inequality to obtain,
uniformly for real~$t$, 
$${1\over N}\sum_{\di{n\leqslant N}{g(n)\leqslant \lambda w+t\sqrt{\gamma w}}}1-\Phi(t)\ll{1\over \sqrt{w}}+\int_{-c_0}^{c_0}\abs{\Delta(\vartheta,w)\over\vartheta}\d\vartheta\ll{1\over \sqrt{\log_2 N}}\cdot$$
However, a standard estimate (see, \eg, \citeplus{Te15}{th. III.3.8})  furnishes that $$g(n)-\lambda w=G^*_N(n;z)+O\big(\log_3N\big)$$ for all but at most $\ll N/\sqrt{\log_2N}$ integers $n\leqslant N$. This plainly establishes \eqref{floiG}.  
\par\medskip
In order to deduce \eqref{evalmoms} from \eqref{floiG}, we need an argument enabling us to discard very large values of $G^*(n;z)$.
To this end, we employ the Laplace transform
$$\L(N;\vartheta):=L(N;-i\vartheta)=\sum_{n\leqslant N}\e^{\vartheta g(n)}\qquad (\vartheta\in\r).$$
Arguing as previously, we get, for $\vartheta\ll1$,
$$\L(N;\vartheta)=N\FF(\vartheta,w)+O\Big({N\over (\log N)^{1/T}}\Big),\eqdef{Laxt2}$$
with now $$\FF(\vartheta,w):=\e^{-\vartheta \lambda w}\sum_{k\leqslant w/z}{\lambda^k(\e^{\vartheta}-1)^k(w-kz)^k\over k!}\cdot\eqdef{defFL}$$
\par 
We  aim to show that
$$
\FF(\vartheta,w)\ll \e^{\fs12\lambda\vartheta^2w}\quad (w^{-1/2}\leqslant |\vartheta|\leqslant w^{-1/3}).\eqdef{majF}$$

\par 
 For $\vartheta\geqslant 0$,  this follows trivially from bounding $w-kz$ by $w$ in \eqref{defF}. The complementary case turns out to be more subtle.  Writing $\nu:=\lambda(1-\e^{-\vartheta})$, we have, for $\alpha>0$, $K=w/z$,
$$\e^{-\vartheta\lambda w}\FF(-\vartheta,w)=\sum_{0\leqslant k\leqslant w/z}{(-1)^k\nu^k(w-kz)^k\over k!}={I(\alpha)-I(-K-\dm)\over 2\pi i},$$
where $$I(a):=\int_{a-i\infty}^{a+i\infty}H(s)\d s\quad \big(a\geqslant -K-\dm\big),\qquad H(s):={\Gamma(s)\over \{\nu(w+sz)\}^s}\cdot$$ This stems from the fact that the residue of $\Gamma(s)$ at $s=-k$ is $(-1)^k/k!$ and from a classical  upper bound for $\Gamma(s)$ in vertical strips.
\par\goodbreak
We select $\alpha:=w\nu/(1-z\nu)$, so that $\nu(w+\alpha z)=\alpha\gg\sqrt{w}$. By the complex Stirling's formula (see, e.g., \citeplus{Te15}{th. II.0.12}), we have, for $ s=\alpha+i\tau$,
$$H(s)\ll{\e^{-\alpha}\over \sqrt{\alpha}}\abs{\bigg({1+i\tau/\alpha\over 1+i\tau\nu z/\alpha}\bigg)^s}.$$
Now
$$\eqalign{\re\big\{(\alpha&+i\tau)\log \big(1+i\tau/\alpha\big)\big\}=\dm\alpha\log \big(\alpha^2+\tau^2\big)-\tau\arctan\Big({\tau\over \alpha}\Big)\cr&\leqslant \dm\alpha\log \Big(1+{\tau^2\over \alpha^2}\Big)-\tau\arctan\Big({\tau\over \alpha}\Big)=\int_0^{|\tau|/\alpha}{\alpha t-|\tau|\over 1+t^2}\d t,\cr}$$
and similarly
$$-\re\big\{(\alpha+i\tau)\log \big( 1+i\tau\nu z/\alpha\big)\big\}\leqslant \int_0^{|\tau|\nu z/\alpha}{|\tau|-\alpha t\over 1+t^2}\d t.$$
Therefore $H(s)\ll\e^{-\alpha}/\sqrt{\alpha}$, and we infer that the contribution to $I(\alpha)$ of the range $|\tau|\leqslant w^2$~is 
$$\ll \e^{-\alpha}w^2/\sqrt{\alpha}\ll\e^{-\alpha}w^{2}.\eqdef{Ialpha1}$$
However, in the case $|\tau|>w^2$, we have, by a further appeal to Stirling's complex formula,
$$H(s)=\sqrt{2\pi\over s}(\nu z\e)^{-s}\Big(1+{w\over zs}\Big)^{-s}\Big\{1+O\Big({1\over s}\Big)\Big\}=\sqrt{2\pi\over s}(\nu z\e)^{-s}\e^{-w/z}\Big\{1+O\Big({w^2\over s}\Big)\Big\}.$$
The main term may be estimated by partial integration, while the remainder furnishes an absolutely convergent contribution. The overall upper bound is $\ll\e^{-w/z+\alpha\log (1/\nu z\e)}$, which is negligible in front of \eqref{Ialpha1}. \par 
Thus we have proved that
$$I(\alpha)\ll \e^{-w \nu/(1-z\nu)}w^2.\eqdef{Ialpha}$$
\par 
To bound $I(-K-\dm)$, we observe that, for $\re s=-\sigma=-K-\dm$ and suitable $c=c_z>0$, we have
$$\abs{\Gamma(s)\over \nu^s(w+sz)^s}=\abs{\nu^\sigma\Gamma(s+K+1)\over (w+sz)^s\prod_{0\leqslant j\leqslant K}(j+\dm+i\tau)}\ll\nu^K\e^{-c|\tau|}, $$
using Stirling's formula on the line $\re s=\dm$ and the estimates
$$\eqalign{&|{(w+sz)^{-s}}|\ll(z^2+\tau^2)^{\sigma/2}\leqslant (z^2+\tau^2)^{K/2+1/4},\quad\prod_{1\leqslant j\leqslant K}\Big\{{z^2+\tau^2\over (j+\dm)^2+\tau^2}\Big\}^{1/2}\ll\e^{-c|\tau|}.\cr}$$
This furnishes the bound $I(-K-\dm)\ll\nu^K$, which, in the considered range, is negligible in front of that of \eqref{Ialpha}.
\par 
Gathering our estimates, we obtain 
$$\FF(-\vartheta,w)\ll\e^{\vartheta\lambda w-w\nu/(1-z\nu)}w^2\ll\e^{\fs12\vartheta^2\gamma w}w^2\ll\e^{\fs12\lambda\vartheta^2w}\qquad (0\leqslant \vartheta\leqslant w^{-1/3}),$$
which finishes the proof of \eqref{majF}.\par 

Since \eqref{majF} implies that the contribution to the $r$th moment of those $n$ with $|G^*_N(n;z)|>T\sqrt{\log_2N}$ is bounded by the claimed error term provided $T>c_r\sqrt{\log_3N}$ for suitable constant $c_r>0$, we obtain \eqref{evalmoms} as an immediate consequence. The case of absolute moments is similar.\bigskip
\bigskip
\goodbreak
\centerline{\twelvebf References}\bigskip
\eightpoint{
\bibtem{De59} H. Delange, Sur des formules dues \`a Atle Selberg, {\it Bull. Sc. Math. \rm 2$^\circ$ s\'erie
\bf 83}, 101--111.\par 
\bibtem{Er46} P. Erd\H os, Some remarks about additive and multiplicative functions, {\it Bull. Amer. Math. Soc. \bf52} (1946), 527Ð537.\par 
\bibtem{Er59} P. Erd\H os,
Some remarks on prime factors of integers.
{\it Canad. J. Math. \bf11} (1959), 161--167. \par
\bibtem{Ru82} I.Z. Ruzsa, Effective results in probabilistic number theory, in : J. Coquet (ed.),
{\it Th\'eorie \'el\'ementaire et analytique des nombres}, D\'ept. Math. Univ. Valenciennes (1982),
107--130.
\bibtem{So21} E. Sofos, Gaps between prime divisors, preprint, June 1st, 2021, https://arxiv.org/abs/2106.00298.
\bibtem{St86} C. Stein, {\it Approximate computation of expectations}, Institute of Mathematical Statistics Lecture Notes-
Monograph Series, 7, Institute of Mathematical Statistics, Hayward, CA, 1986.\par 
\bibtem{Te93} G. Tenenbaum, Cribler les entiers sans grand facteur premier, in: R. C. Vaughan (ed.), Theory
and applications of numbers without large prime factors, {\it Phil. Trans.  R.  Soc. Lond.} A {\bf
345} (1993), 377--384.\par 
\bibtem{Te15} G. Tenenbaum, {\it Introduction to analytic and probabilistic number theory}, 3rd ed., Graduate Studies in Mathematics 163, Amer. Math. Soc. (2015).\par
\bibtem{Te19} G. Tenenbaum, A note on the normal largest  gap between prime factors, {\it J. ThŽor. Nombres Bordeaux \bf31}, \numero3 (2019), 747-749.
\par }
\smallskip
{\leftskip2mm\rightskip-2cm\sevenrm
\gutter=4cm \doublecolumns
 \obeylines \baselineskip=7pt
RŽgis de la Bretche
UniversitŽ de Paris, Sorbonne UniversitŽ, CNRS, 
Institut de Math. de Jussieu-Paris Rive Gauche,
 F-75013 Paris,  
France
\smallskip
{\seventt regis.de-la-breteche@imj-prg.fr}
 GŽrald Tenenbaum\par
Institut \'Elie Cartan\par 
Universit\'e de Lorraine\par
 BP 70239\par
54506 Vand\oe uvre-ls-Nancy Cedex\par
 France
\smallskip
{\seventt gerald.tenenbaum@univ-lorraine.fr}
\singlecolumn
\par}
\bye